# GENERALIZED STOCHASTIC DIFFERENTIAL UTILITY AND PREFERENCE FOR INFORMATION[1]

By Ali Lazrak

*University of British Columbia and Université d'Evry*

This paper develops, in a Brownian information setting, an approach for analyzing the preference for information, a question that motivates the stochastic differential utility (SDU) due to Duffie and Epstein [*Econometrica* **60** (1992) 353–394]. For a class of backward stochastic differential equations (BSDEs) including the generalized SDU [Lazrak and Quenez *Math. Oper. Res.* **28** (2003) 154–180], we formulate the information neutrality property as an invariance principle when the filtration is coarser (or finer) and characterize it. We also provide concrete examples of heterogeneity in information that illustrate explicitly the nonneutrality property for some GSDUs. Our results suggest that, within the GSDUs class of intertemporal utilities, risk aversion or ambiguity aversion are inflexibly linked to the preference for information.

**1. Introduction.** The study of decision making is fundamental to many applications in economics and finance. The decision maker typically faces uncertainty about results of an experiment such as the profitability of a new product or a financial strategy, efficacy of a monetary policy or a social program, state of health and so on. Since many decades, decision theorists have developed theories and tools which help us to think about decision under uncertainty. The ultimate objective of this line of literature in social science is to provide explanations of the behavior under uncertainty and to give a *rational* support for the observable behavior in various contexts.

This paper studies the preference for information for a specific class of intertemporal utilities. For a fixed consumption horizon $T > 0$, a utility

Received August 2002; revised January 2004.
[1]Supported by Université d'Evry, MITACS and the Sauder School of Buisiness at the University of British Columbia and the Natural Sciences and Engineering Research Council of Canada.
*AMS 2000 subject classifications.* 60H10, 60H30.
*Key words and phrases.* Generalized stochastic differential utility, Brownian filtration, information, backward stochastic differential equation.

function is a function mapping the set of objects of choice, that is the pairs of state contingent consumption process $c = \{c_t, 0 \leq t \leq T\}$ and information filtration $\mathcal{A} = \{\mathcal{A}_t, 0 \leq t \leq T\}$ satisfying *the usual conditions*, into $\mathbb{R}$. The question of preference for information consists of analyzing the dependency of a utility function in its filtration information argument. This specific question has been greatly simplified within the familiar context of the von Neumann–Morgenstern expected utility function. Specifically, an expected utility function is defined by

$$(1) \qquad U_t^{\mathcal{A}}(c) = E\left[\int_t^T e^{-\beta(s-t)} v(c_s)\, ds \Big| \mathcal{A}_t \right],$$

for time $t < T$, $v(\cdot)$ is the felicity function and the expectation $E$ is conditioned by the time $t$ available information $\mathcal{A}_t$. In fact, the expected utility model (1) does not allow for preference for information in the sense that if $\mathcal{A} \subsetneq \mathcal{B}$ are two filtrations such that both $\mathcal{A}_0$ and $\mathcal{B}_0$ are trivial, then $U_0^{\mathcal{A}}(c) = U_0^{\mathcal{B}}(c)$ for any $\mathcal{A}$-adapted consumption process $c$.

However, it is often observed that preference for information is relevant in various decision making situations. For instance, in many medical decisions (choice between various form of prenatal diagnosis such as the amniocentesis or decision to test for diseases such as multiple sclerosis), the decision maker must decide whether she wishes to have the true state of health revealed earlier or later. More generally, psychologists have recognized the importance of the feelings related to the prospect of information acquisition [see Grant, Kajii and Polack (1998) and Chew and Ho (1994) and the references cited therein]. It has been recognized that information acquisition has an *extrinsic* and an *intrinsic* motivation.

The extrinsic motivation corresponds to the notion that people value information to take appropriate contingents decisions and thus influence in a favorable way the final outcome. For example, certain medical treatments may lower the severity of a disease which provide an incentive to gather information about the health state. In an investment context for instance, information enhances the planning and should help to identify financial strategies which provide a higher expected profitability. In particular, for an expected utility maximizer, Epstein (1980) has shown in an investment problem that the prospect of greater future information increases the incentives to maintain some flexibility in order to take advantage of the content of the future information. [Yet, as noted above and clarified below, the expected utility investors of Epstein (1980) are indifferent to information and are interested in it only for its planning benefits.]

On the other hand, intrinsic motivation corresponds to the notion that, planning benefits notwithstanding, people like (or dislike) information for its own sake. In other words, intrinsic attitude toward information is defined as individual's direct interest to have access to more (or less) information because they perceive it to be innately satisfying (or unsatisfying). For instance, a decision maker could be intrinsically information lover because she is anxious and prefer to know earlier the outcome of any uncertainty (think of the example of a pregnant women who decide to do the amniocentesis prenatal diagnosis). Another decision maker may be intrinsically information averse since he fears a bad outcome (think to the example of a person who delay a test for disease) or simply because he is optimistic and prefers a hopeful feeling rather than risking a sad news.

In a pioneering discrete time model, Kreps and Porteus (1978), generalized the von Neumann–Morgenstern expected utility model (1) to permit intrinsic information aversion or information loving. This work gave rise to the stochastic differential utility (SDU) [Duffie and Epstein (1992) and its discrete time counterpart Epstein and Zin (1989)]. The SDU generalized the expected utility model (1) and is associated with an "intertemporal aggregator" $f$, a function satisfying appropriate conditions. The SDU is defined by

$$(2) \qquad U_t^{\mathcal{A}}(c) = E\left[\int_t^T f(c_s, U_s^{\mathcal{A}}(c))\,ds \Big| \mathcal{A}_t\right].$$

The SDU model reduces to the additive model (1) when $f$ is linear, $f(c,u) = v(c) - \beta u$. The SDU was primarily motivated by the desire to have some flexibility in the modeling of the concepts of risk aversion and the concept of consumption intertemporal substitution. While the two concepts were governed by the same parameter in the time additive model (1), the SDU allowed some separation of these two aspects of the preferences. This feature was particularly relevant from an empirical perspective since it helped to match more closely consumption rates data and equity returns data in the US [Epstein and Zin (1991)]. At the same time, unlike the expected utility model (1), the SDU model (2) exhibits an intrinsic attitude toward information. From a mathematical perspective, intrinsic attitude toward information is characterized by the fact that the initial value of the SDU (2) [i.e., $U_0^{\mathcal{A}}(c)$] depends not only on consumption but also on the filtration $\mathcal{A}$.

Building on the discrete time approach to intrinsic attitude toward information of Kreps and Porteus (1978), Skiadas (1998) shows that in the continuous time SDU model (2), the concavity (convexity) of an intertemporal aggregator with respect to its utility argument $U$ implies an intrinsic preference for late (early) resolution of uncertainty. To illustrate their point, consider two filtration $\mathcal{F} \subset \mathcal{G}$ and an intertemporal aggregator $f$ which is concave with respect to its utility argument $U$. Then, Jensen's inequality gives,

$$E(U_t^{\mathcal{G}}(c)|\mathcal{F}_t) = E\left[\int_t^T f(c_s, U_s^{\mathcal{G}}(c))\,ds \Big| \mathcal{F}_t\right]$$

$$\leq E\left[\int_t^T f(c_s, E(U_s^{\mathcal{G}}(c)|\mathcal{F}_s))\,ds\bigg|\mathcal{F}_t\right],$$

for any consumption process $c$ which is progressively measurable with respect to the coarser filtration $\mathcal{F}$. Thus the optional projection process $E[U_\cdot^{\mathcal{G}}(c)|\mathcal{F}_\cdot]$ maybe interpreted as a sub-solution of the recursion (2) in the setting of the filtration $\mathcal{A} = \mathcal{F}$ and as such, heuristically, the sub-solution $E[U_\cdot^{\mathcal{G}}(c)|\mathcal{F}_\cdot]$ is smaller ($P \otimes dt$ a.s.) than the solution itself $U_\cdot^{\mathcal{F}}(c)$. Unfortunately, while pointing an elegant way of proving out an elegant way of proving a monotonicity of a utility functional with respect to its filtration argument, the Kreps–Porteus–Skiadas method only provides sufficient conditions for preference for early (or late) resolution of uncertainty and no characterization is obtained.

This paper is an attempt to analyze the question of intrinsic attitude toward information within a more general class of utility functions. We consider the class of generalized stochastic differential utility (GSDU) introduced in Lazrak and Quenez (2003). It has been shown by Skiadas (2003) and Lazrak and Quenez (2003) that the GSDU unifies the SDU of Duffie and Epstein (1992) and a recent class of intertemporal utility functions. This class encompasses the portfolio decision models of Chen and Epstein (2002) and of Anderson, Hansen and Sargent (1998). These models have been introduced with the operational objective of modeling the imperfect knowledge of the asset returns probability distribution and its impact on portfolio decision and asset prices. The objective of this paper is to identify the implicit implications of these utilities from the angle of the intrinsic attitude toward information. More specifically, we took the view that investors have a neutral intrinsic attitude toward information (in a sense to be made precise later) and, characterize this property in a context of Brownian information and under certain assumptions on the (generalized) intertemporal aggregator. Our finding suggests that, in general, the information neurality will not hold for our class of GSDU. However, when information heterogeneity is such that the Brownian property is preserved under the finer filtration, neutrality for information holds.

Therefore, the GSDUs class of utility functions are generally not information neutral and this suggests that the risk attitudes and the ambiguity attitudes are in some sense confounded with the information attitude within this class of utility functions. Consequently, any prediction of these models for portfolio decision or asset prices is also induced by the extent to which these utilities exhibit preference for information. Finally, our results should be of interest to the literature on the design of risk measure for institutional investors and financial institutions [see Artzner, Delbaen, Eber and Heath (1999), Wang (2000), Artzner, Delbaen, Eber, Heath and Ku (2002) and Riedel (2002)]. In fact, a GSDU is in some sense a dynamic risk measure and a preference for information may be desirable in that context. For instance, in a stock portfolio management context, it is possible to have a view about how a risk measure should depend on the timing of information releases on the stock prices. In particular, the GSDU would have then the ability to provide a quantitative prediction of the utility cost of an information enhancement such as an increase of the frequency of the accounting reports of the underlying companies or perhaps an increase of the coverage of the financial analysts.

The paper is organized as follows. In Section 2 we give the exact setting for two generalized versions of model (2), define the information neutrality property and give some mathematical prerequisites. Section 3 gives some concrete filtrations and utility models encompassed by our formulation. Section 4 develops some GSDU computations for two examples of heterogeneity in information that illustrate the problem. The first example (Brownian anticipation) exhibits a situation where information neutrality does not hold. In the second example, the coarser filtration is generated by the absolute value of the Brownian motion that drives the finer filtration and we will see that in this context the information neutrality may hold. In Section 5, we characterize the information neutrality for a class of BSDEs (including GSDUs) which driver depend on intensity $Z$ only through its Euclidean norm $\|Z\|$. We show that the second example of Section 4 offers the only type of information heterogeneity that allows information neutrality to hold for this class of BSDEs. In Section 6, we conclude.

## 2. The model.

2.1. *Context and definitions.* Let $(\Omega, \mathcal{G}, P)$ be a complete probability space and for the fixed time $T$, let $\mathcal{G}_{(\cdot)} = \{\mathcal{G}_t, 0 \leq t \leq T\}$ and $\mathcal{F}_{(\cdot)} = \{\mathcal{F}_t, 0 \leq t \leq T\}$ be two filtrations that contain all negligible events and are right-continuous and such that $\mathcal{F}_{(\cdot)} \subsetneq \mathcal{G}_{(\cdot)}$. Furthermore, we suppose that the filtration $\mathcal{G}_{(\cdot)}$ (resp. $\mathcal{F}_{(\cdot)}$) has a predictable representation property with respect to a standard $n$-dimensional Brownian motions $W^{\mathcal{G}} = (^1W^{\mathcal{G}}, ^2W^{\mathcal{G}}, \ldots, ^nW^{\mathcal{G}})$ [resp. $W^{\mathcal{F}} = (^1W^{\mathcal{F}}, ^2W^{\mathcal{F}}, \ldots, ^nW^{\mathcal{F}})$]: For $\mathcal{A} \in \{\mathcal{G}, \mathcal{F}\}$, each $\mathcal{A}$-local martingale $M$ can be represented as a stochastic integral with respect to $W^{\mathcal{A}}$, that is, there exists an $\mathcal{A}$-predictable process $\varphi$ in $\mathbb{R}^n$ with $\int_0^T \|\varphi\|^2 \, dt < \infty$ a.s. such that $M_t = M_0 + \int_0^t \varphi_s \cdot dW_s^{\mathcal{A}}$, $0 \leq t \leq T$. In other words, following the Revuz and Yor terminology [e.g., Revuz and Yor (1999), page 219], the filtrations $\mathcal{G}$ and $\mathcal{F}$ are weakly Brownian. As we will illustrate with some specific examples in Section 3, there are many ways of constructing such a couple of filtrations (representing heterogeneous information). It is important to notice at this stage that in general, the process $W^{\mathcal{F}}$ is not a Brownian motion under the finer filtration $\mathcal{G}$. However, as we shall illustrate in Section 3, there are some special specifications of the filtrations $\mathcal{G}$ and $\mathcal{F}$ under which the process $W^{\mathcal{F}}$ turns out to be a $\mathcal{G}$-Brownian motion.

We shall denote by $\mathcal{P}^{\mathcal{G}}$ the $\mathcal{G}_{(\cdot)}$-predictable $\sigma$-field and by $\mathcal{P}^{\mathcal{F}}$ the $\mathcal{F}_{(\cdot)}$-predictable $\sigma$-field. We consider for each integer $p$ the sets $\mathcal{H}^2(\mathcal{G}, \mathbb{R}^p) = \{X : [0, T] \times \Omega \to \mathbb{R}^p / X \in \mathcal{P}^{\mathcal{G}}$ and $E[\int_0^T |X_s|^2 \, ds] < \infty\}$ and $\mathcal{H}^2(\mathcal{F}, \mathbb{R}^p) = \mathcal{H}^2(\mathcal{G}, \mathbb{R}^p) \cap \mathcal{P}^{\mathcal{F}}$.

For each random variable $\xi \in L^2(\mathcal{F}_T)$, we define the BSDE $Y^{\mathcal{G}}(\xi) \in \mathcal{H}^2(\mathcal{G}, \mathbb{R})$ associated to the filtration $\mathcal{G}_{(\cdot)}$ as the solution of the recursion

$$
\begin{aligned}
Y_t^{\mathcal{G}}(\xi) &= \xi + \int_t^T h(s, \omega, Y_s^{\mathcal{G}}(\xi), Z_s^{\mathcal{G}}(\xi)) \, ds - \int_t^T Z_s^{\mathcal{G}}(\xi) \cdot dW_s^{\mathcal{G}} \\
&\equiv E\bigg[\xi + \int_t^T h(s, \omega, Y_s^{\mathcal{G}}(\xi), Z_s^{\mathcal{G}}(\xi)) \, ds \bigg| \mathcal{G}_t \bigg],
\end{aligned}
\tag{3}
$$

where *the driver* $h$ defined on $[0, T] \times \Omega \times \mathbb{R} \times \mathbb{R}^n$ with values in $\mathbb{R}$, s.t. $(h(t, \omega, y, z))_{0 \leq t \leq T} \in \mathcal{H}^2(\mathcal{F}, \mathbb{R})$ for each $(y, z) \in \mathbb{R} \times \mathbb{R}^n$ and $h$ satisfies the following standing assumptions.

*Standing assumptions.* (A1) There exists a constant $K \geq 0$ s.t. $P$-a.s., we have

$$\forall t, \ \forall (y_1, y_2), \ \forall (z_1, z_2)$$
$$|h(s, \omega, y_1, z_1) - h(s, \omega, y_2, z_2)| \leq K(|y_1 - y_2| + \|z_1 - z_2\|).$$

(A2) The process $(h(t, \omega, 0, 0))_{0 \leq t \leq T}$ belongs to $\mathcal{H}^2(\mathcal{F}, \mathbb{R})$.

Note that the process $Z^{\mathcal{G}}(\xi) \in \mathcal{H}^2(\mathcal{G}, \mathbb{R}^n)$ is part of the solution of (3) and we call it the intensity associated to the BSDE (3).

Similarly, we define the BSDE $Y^{\mathcal{F}}(\xi) \in \mathcal{H}^2(\mathcal{F}, \mathbb{R})$ associated to the filtration $\mathcal{F}_{(\cdot)}$ as the solution of the BSDE

$$
\begin{aligned}
Y_t^{\mathcal{F}}(\xi) &= \xi + \int_t^T h(s, \omega, Y_s^{\mathcal{F}}(\xi), Z_s^{\mathcal{F}}(\xi)) \, ds - \int_t^T Z_s^{\mathcal{F}}(\xi) \cdot dW_s^{\mathcal{F}} \\
&\equiv E\bigg[\xi + \int_t^T h(s, \omega, Y_s^{\mathcal{F}}(\xi), Z_s^{\mathcal{F}}(\xi)) \, ds \bigg| \mathcal{F}_t \bigg].
\end{aligned}
\tag{4}
$$

We will also be interested by a second class of BSDEs, GSDU, an extension of Duffie and Epstein (1992) model of SDU that has been proposed in Lazrak and Quenez (2003). For any given *contingent consumption plan*, a process $c \in \mathcal{H}^2(\mathcal{F}, \mathbb{R})$, the GSDU $U^{\mathcal{G}}(c) \in \mathcal{H}^2(\mathcal{G}, \mathbb{R})$ associated to the filtration $\mathcal{G}_{(\cdot)}$ solves the recursion

$$
\begin{aligned}
U_t^{\mathcal{G}}(c) &= \int_t^T f(s, c_s, U_s^{\mathcal{G}}(c), V_s^{\mathcal{G}}(c)) \, ds - \int_t^T V_s^{\mathcal{G}}(c) \cdot dW_s^{\mathcal{G}} \\
&\equiv E\bigg[\int_t^T f(s, c_s, U_s^{\mathcal{G}}(c), V_s^{\mathcal{G}}(c)) \, ds \bigg| \mathcal{G}_t \bigg],
\end{aligned}
\tag{5}
$$

where *the intertemporal aggregator* $f$ defined on $[0,T] \times \mathbb{R} \times \mathbb{R} \times \mathbb{R}^n$ with values in $\mathbb{R}$, s.t. $(f(t, c_t, y, z))_{0 \leq t \leq T} \in \mathcal{H}^2(\mathcal{F}, \mathbb{R})$ for each $(y, z) \in \mathbb{R} \times \mathbb{R}^n$ and $f$ satisfies the following standing assumptions:

*Standing assumptions.* (B1) There exists a constant $K \geq 0$ s.t., $P$-a.s., for all relevant $(t, c, y_1, y_2, z_1, z_2)$ we have

$$|f(t, c, y_1, z_1) - f(t, c, y_2, z_2)| \leq K(|y_1 - y_2| + \|z_1 - z_2\|).$$

(B2) There exists some positive constants $k_1, k_2$ and $0 < p < 1$ s.t. $|f(t, c, 0, 0)| \leq k_1 + k_2 c^p$.

In fact, Duffie and Epstein (1992) define SDU of the form (5) in a context where the intertemporal aggregator is essentially independent of $z$, and thus we shall call this case the *classical SDU*.

We define as well the GSDU associated to the filtration $\mathcal{F}_{(\cdot)}$ as the solution of the BSDE

$$\begin{aligned}
U_t^{\mathcal{F}}(c) &= \int_t^T f(s, c, U_s^{\mathcal{F}}(c), V_s^{\mathcal{F}}(c))\, ds - \int_t^T V_s^{\mathcal{F}}(c) \cdot dW_s^{\mathcal{F}} \\
&\equiv E\left[\int_t^T f(s, c_s, U_s^{\mathcal{F}}(c), V_s^{\mathcal{F}}(c))\, ds \bigg| \mathcal{F}_t\right]
\end{aligned}$$
(6)

for each $c \in \mathcal{H}^2(\mathcal{F}, \mathbb{R})$.

Note that the intertemporal aggregator $f$ is a deterministic function of $(t, c, y, z)$ and thus the GSDU model (5) is a special case of the BSDE model (3) that is obtained formally by setting $\xi = 0$ and $h(t, \omega, y, z) = f(t, c_t, y, z)$. However, as will be seen from the following definitions, the information neutrality property has a different meaning in the two models and therefore a different method is needed to characterize it in the two models.

Now let us define the information neutrality property.

DEFINITION 1. A BSDE exhibits *information neutrality* (between the filtration $\mathcal{F}$ and $\mathcal{G}$) if and only if the solutions of the BSDEs (3) and (4) satisfy

(7) $$Y_t^{\mathcal{F}}(\xi) = Y_t^{\mathcal{G}}(\xi), \qquad dP \otimes dt \text{ a.s.,}$$

for all $\xi \in L^2(\mathcal{F}_T)$.

DEFINITION 2. A GSDU exhibits *information neutrality* (between the filtration $\mathcal{F}$ and $\mathcal{G}$) if and only if the solutions of the GSDUs (5) and (6) satisfy

(8) $$U_t^{\mathcal{F}}(c) = U_t^{\mathcal{G}}(c), \qquad P \otimes dt \text{ a.s.,}$$

for all $c \in \mathcal{H}^2(\mathcal{F}, \mathbb{R})$.

We motivate these definitions by interpreting (7) and (8) as expressing an indifference for the purpose of decision toward the otherwise anticipated utility of more rather less information. A decision maker who exhibits such a property has no intrinsic motivation to gather information for a fixed consumption. In the subsequent analysis, our objective is to characterize this property.

2.2. *Mathematical background.* Under our Lipschitz assumptions on the driver/aggregator [assumptions (A1) and (B1)], it is now a standard result in the BSDE literature that existence and uniqueness (in a suitable sense) hold for the recursions (3)–(6).

More precisely, for $\mathcal{A} \in \{\mathcal{G}, \mathcal{F}\}$, it follows from Pardoux and Peng (1990) [see also El Karoui, Peng and Quenez (1997) and Ma and Yhong (1999)] that under assumptions (A1) and (A2) and, for each $\xi \in L^2(\mathcal{F}_T)$, there exist a unique pair $(Y^{\mathcal{A}}(\xi), Z^{\mathcal{A}}(\xi)) \in \mathcal{H}^2(\mathcal{A}, \mathbb{R}) \times \mathcal{H}^2(\mathcal{A}, \mathbb{R}^n)$ such that

$$(9) \qquad Y_t^{\mathcal{A}}(\xi) = E\bigg[\xi + \int_t^T h(s, \omega, Y_s^{\mathcal{A}}(\xi), Z_s^{\mathcal{A}}(\xi))\, ds \bigg| \mathcal{A}_t\bigg].$$

Similarly, for $\mathcal{A} \in \{\mathcal{G}, \mathcal{F}\}$, it follows from Pardoux and Peng (1990) that under assumptions (B1) and (B2) and, for each $c \in \mathcal{H}^2(\mathcal{F}, \mathbb{R})$ there exist a unique pair $(U^{\mathcal{A}}(c), V^{\mathcal{A}}(c)) \in \mathcal{H}^2(\mathcal{A}, \mathbb{R}) \times \mathcal{H}^2(\mathcal{A}, \mathbb{R}^n)$ such that

$$(10) \qquad U_t^{\mathcal{A}}(c) = E\bigg[\int_t^T f(s, c_s, U_s^{\mathcal{A}}(c), V_s^{\mathcal{A}}(c))\, ds \bigg| \mathcal{A}_t\bigg].$$

Now, since we will extensively use them in the subsequent analysis, it is worthwhile to recall [see El Karoui, Peng and Quenez (1997)] the representation theorems of linear (resp. concave) BSDEs as a conditional expectation (resp. an essential infinimum of conditional expectations). We will state these results for a filtration $\mathcal{A} \in \{\mathcal{G}, \mathcal{F}\}$ and only for the BSDE model (9) [the GSDU model (10) being a particular case of the BSDE model (9)].

PROPOSITION 1. *Let $(\rho, \kappa)$ be a bounded $(\mathbb{R}, \mathbb{R}^n)$-valued $\mathcal{A}$-predictable process, $\varphi$ an element of $\mathcal{H}^2(\mathcal{A}, \mathbb{R})$ and $\xi$ and element of $L^2(\mathcal{A}_T)$. Then the linear BSDE*

$$Y_t^{\mathcal{A}}(\xi) = \xi + \int_t^T (\varphi_s + Y_s^{\mathcal{A}}(\xi)\rho_s + Z_s^{\mathcal{A}}(\xi) \cdot \kappa_s)\, ds - \int_t^T Z_s^{\mathcal{A}}(\xi) \cdot dW_s^{\mathcal{A}}$$

*has a unique solution $(Y^{\mathcal{A}}(\xi), Z^{\mathcal{A}}(\xi)) \in \mathcal{H}^2(\mathcal{A}, \mathbb{R}) \times \mathcal{H}^2(\mathcal{A}, \mathbb{R}^n)$ which admits the representation*

$$Y_t^{\mathcal{A}}(\xi) = E\bigg[\Upsilon_t^T \xi + \int_t^T \Upsilon_t^s \varphi_s\, ds \bigg| \mathcal{A}_t\bigg],$$

*where $\Upsilon_t^s$ is the adjoint process defined for $s \geq t$ by the forward SDE*

$$d\Upsilon_t^s = \Upsilon_t^s[\rho_s\, ds + \kappa_s \cdot dW_s^{\mathcal{A}}], \qquad \Upsilon_t^t = 1.$$

Alternatively, when the driver $h$ is concave with respect to $(y, z)$, it is possible to express it as an infinimun of linear functions of $(y, z)$: denoting by $H$ the polar function of $h$ defined by

$$H(t, \rho, \kappa) = \sup_{(y,z) \in \mathbb{R} \times \mathbb{R}^n} [h(t, y, z) - \rho y - \kappa \cdot z],$$

the conjugacy relationship gives [for each $(\omega, t)$]

$$h(t, y, z) = \inf_{(\rho, \kappa) \in [-K, K]^{n+1}} [h^{(\rho, \kappa)}(t, y, z)],$$

where $h^{(\rho,\kappa)}(t, y, z) = H(t, \rho, \kappa) - \rho y - \kappa \cdot z$ and where we recall that $K$ is the Lipschitz constant for the driver $h$ (the domain of definition of $H$ is a subset of $[-K, K]^{n+1}$). Heuristically, the representation theorem for concave BSDEs states that the infinimum of the above conjugacy relationship commutes with the BSDE transform, that is,

$$Y_t(h) \equiv Y_t(\operatorname{ess\,inf} h^{(\rho,\kappa)}) = \operatorname{ess\,inf} Y_t(h^{(\rho,\kappa)}).$$

In order to state this result more precisely in the following proposition, we first define the domain [see El Karoui, Peng and Quenez (1997)]

$$\mathcal{D} := \{(\rho, \kappa) \in \mathcal{P}^{\mathcal{A}} \cap [-K, K]^{n+1} | H(\cdot, \rho_\cdot, \kappa_\cdot) \in \mathcal{H}^2(\mathcal{A}, \mathbb{R})\}.$$

PROPOSITION 2. *Let $h$ be a concave driver satisfying assumptions* (A1) *and* (A2) *and let $H$ the associated polar function. Then the BSDE* (9) *admits the dual representation*

$$Y_t^{\mathcal{A}}(\xi) = \operatorname{ess\,inf}_{(\rho, \kappa) \in \mathcal{D}} E\left[\Upsilon_{t,T}^{\rho, \kappa} \xi + \int_t^T \Upsilon_{t,s}^{\rho, \kappa} H(s, \rho_s, \kappa_s) \, ds \Big| \mathcal{A}_t\right],$$

*where $\Upsilon_{t,s}^{\rho,\kappa}$ is the adjoint process defined for $s \geq t$ by the forward SDE*

$$d\Upsilon_{t,s}^{\rho,\kappa} = \Upsilon_{t,s}^{\rho,\kappa} [\rho_s \, ds + \kappa_s \cdot dW_s^{\mathcal{A}}], \qquad \Upsilon_{t,t}^{\rho,\kappa} = 1.$$

## 3. Some examples of filtrations and utilities.

3.1. *Examples of heterogeneous filtrations.* There are *many* ways to construct a sequence of coarser or finer Brownian filtrations and we give here some examples of constructions.

*Losing the sign of a Brownian motion.* Departing from a completed filtration $\mathcal{B}$ generated a two-dimensional Brownian motion $(\nu_t^{\mathcal{B}} = (\nu_{1t}^{\mathcal{B}}, \nu_{2t}^{\mathcal{B}}))_{0 \leq t \leq T}$ one can construct the filtration $\mathcal{A}$ generated by

$$\mathcal{A}_t := \sigma((|\nu_{1s}^{\mathcal{B}}|, |\nu_{2s}^{\mathcal{B}}|); 0 \leq s \leq t),$$

and it follows from Revuz and Yor (1999) that $\mathcal{A}$ is generated by the two dimensional Brownian motion $(\nu_t^{\mathcal{A}})_{0 \leq t \leq T} \equiv (\int_0^t \operatorname{sgn}(\nu_{1s}^{\mathcal{B}}) \, d\nu_{1s}^{\mathcal{B}}, \int_0^t \operatorname{sgn}(\nu_{2s}^{\mathcal{B}}) \, d\nu_{2s}^{\mathcal{B}})_{0 \leq t \leq T}$.

Note that this method provides a way to construct an infinite sequence of coarser filtrations. Finally, it is important to observe that in this particular example, $\nu^{\mathcal{A}}$ is also a Brownian motion under the filtration $\mathcal{B}$.

*Brownian anticipation.* Consider an $n$-dimensional Brownian motion $(\nu_s^{\mathcal{A}}; 0 \leq s \leq t)$ that generates a completed filtration $\mathcal{A}$. Then the process

$$\nu_t^{\mathcal{B}} := \tfrac{1}{\sqrt{2}} \nu_{2t}^{\mathcal{A}}, \qquad 0 \leq t \leq T,$$

is a Brownian motion and generates a completed filtration $\mathcal{B}$ that satisfies

$$\mathcal{B}_t = \mathcal{A}_{2t} \supset \mathcal{A}_t.$$

Notice that in this context, $\nu$ is not a $\mathcal{B}$-martingale and thus it is not a Brownian motion under $\mathcal{B}$.

*Brownian motion with an independent random drift.* A third example comes from filtering theory. Consider a scalar Brownian motion $(\nu_t)_{0 \leq t \leq T}$ and an independent and integrable random variable $\mu$. Consider the filtration $\mathcal{B}$ generated by

$$\mathcal{B}_t := \sigma(\nu_s; 0 \leq s \leq t) \vee \sigma(\mu)$$

and its subfiltration

$$\mathcal{A}_t := \sigma(\nu_s + \mu s; 0 \leq s \leq t),$$

which is well known to be generated [see, e.g., Liptser and Shiryayev (1977)] by the Brownian motion $(\nu_t^{\mathcal{A}})_{0 \leq t \leq T} = (\nu_t + \int_0^t (\mu - E(\mu|\mathcal{A}_s))\,ds)_{0 \leq t \leq T}$. Note however that this example is outside the scope of this paper since $\mathcal{B}_0$ is not trivial.

3.2. *Examples of intertemporal aggregator.*

*The multi-prior expected utility process.* When the BSDE driver has the form

(11) $$h(t, \omega, y, z) = -\sum_{i=1}^n k_i |z_i|,$$

where $k_i > 0$ for $i = 1, \ldots, n$, the process $Y^{\mathcal{F}}(U(\xi))$ defined in (3) [resp. the process $Y^{\mathcal{G}}(U(\xi))$ defined in (4)] with the terminal data $U(\xi)$ where $U(\cdot)$ is a nondecreasing and concave function mapping $\mathbb{R}$ onto $\mathbb{R}$ may be interpreted as a multi-prior utility for the wealth $\xi$ [Chen and Epstein (2002)].

Alternatively, in the GSDU case, when the intertemporal aggregator has the following form:

(12) $$f(t, c, y, z) = u(c) - \sum_{i=1}^n k_i |z_i|,$$

where $k_i > 0$ for $i = 1, \ldots, n$ and where $u(\cdot)$ is a nondecreasing and concave function mapping $\mathbb{R}$ onto $\mathbb{R}$, the process $U^{\mathcal{F}}(c)$ defined in (5) [resp. the process $U^{\mathcal{G}}(c)$ defined in (6)] may be interpreted as a multi-prior utility for the consumption process $c$ [Chen and Epstein (2002)].

*Quadratic GSDU.* When the intertemporal aggregator has the form

$$f(c, y, z) = \log(c) - \beta y - \frac{\alpha}{2} z^2, \tag{13}$$

with parameter restrictions: $\beta \geq 0$ and $\alpha \geq 0$ the existence of the GSDUs (5) are not guaranteed anymore since the intertemporal aggregator $f$ is not Lipschitz with respect to $z$. The GSDU associated to (13) is in fact in the class of quadratic BSDE that has been extensively studied in Kobylansky (2000) who shows the existence by an approximation technique. In the specific case under consideration, Schroder and Skiadas (1999) show the existence and uniqueness of the solution of the BSDE (5) and consequently of the BSDE (6). Their proof consists of building an appropriate set of consumption plan processes (that contains but is not limited to the set of bounded process) and involves a fixed point theorem.

Interestingly, the model (13) is important since it has been recently shown by Skiadas (2003) and Lazrak and Quenez (2003) that the GSDU associated with (13) is a unified formulation of a recent approach to uncertainty aversion related to the robust control theory. This approach has been introduced by Anderson, Hansen and Sargent (1998) [see also Hansen, Sargent, Turmuhambetova and Williams (2002) and Uppal and Wang (2003) for some applications of that model to asset pricing issues].

**4. An illustrative example.** The objective of this section is to give a concrete situation where we can measure explicitly the utility under heterogeneous filtrations. In fact, the example that we shall give in the sequel was very helpful to us as a guide of how to handle the problem given in the previous section.

Assume that the filtration $\mathcal{F}$ is generated by

$$\mathcal{F}_t := \sigma(|W_s^{\mathcal{G}}|; 0 \leq s \leq t) = \sigma(W_s^{\mathcal{F}}; 0 \leq s \leq t), \tag{14}$$

where the $\mathcal{F}$-Brownian motion $(W_t^{\mathcal{F}}; 0 \leq t \leq T)$ is given by

$$W_t^{\mathcal{F}} = \int_0^t \operatorname{sgn}(W_s^{\mathcal{G}}) \, dW_s^{\mathcal{G}}.$$

Also, we will consider the anticipating filtration

$$\mathcal{H}_t := \sigma(W_{2s}^{\mathcal{F}}; 0 \leq s \leq t) = \sigma(W_s^{\mathcal{H}}; 0 \leq s \leq t),$$

where the $\mathcal{H}$-Brownian motion $(W_t^{\mathcal{H}}; 0 \leq t \leq T)$ is given by

$$W_t^{\mathcal{H}} = \frac{1}{\sqrt{2}} W_{2t}^{\mathcal{F}}.$$

It is clear that $\mathcal{F}_{(\cdot)} \subsetneq \mathcal{G}_{(\cdot)}$ and $\mathcal{F}_{(\cdot)} \subsetneq \mathcal{H}_{(\cdot)}$, and in order to simplify the exposition, assume furthermore that the Brownian motion $W^{\mathcal{G}}$ is one dimensional ($n = 1$). Now, in the following sections, we shall consider the particular consumption plan $b \in \mathcal{H}^2(\mathcal{F}, \mathbb{R})$ given by

$$b_t = \exp(W_t^{\mathcal{F}}), \tag{15}$$

and compute its associated GSDU for some simple intertemporal aggregators under the filtrations $\mathcal{F}_{(\cdot)}, \mathcal{G}_{(\cdot)}$ and $\mathcal{H}_{(\cdot)}$.

4.1. *A linear GSDU intertemporal aggregator.* Now, let us analyze what happens if we introduce the *simplest* dependence in $z$ in the intertemporal aggregator, that is, a linear additive dependence of the form

$$(16) \qquad f(s,c,y,z) = \log(c) - \beta(s)y - \gamma z,$$

where the $\beta(\cdot)$ is a deterministic integrable function.

From Proposition 1 and for $\mathcal{A} \in \{\mathcal{F}, \mathcal{G}, \mathcal{H}\}$, the GSDU associated with the aggregator (16) admits the representation

$$(17) \qquad U_t^{\mathcal{A}}(b) = E\bigg[\int_t^T \Upsilon_t^s W_s^{\mathcal{F}} \, ds \bigg| \mathcal{A}_t\bigg],$$

where $\Upsilon_t^s$ is the adjoint process defined for $s \geq t$ by the forward SDE

$$d\Upsilon_t^s = -\Upsilon_t^s[\beta(s)\,ds + \gamma\,dW_s^{\mathcal{A}}], \qquad \Upsilon_t^t = 1.$$

When $\gamma = 0$, the adjoint process $\Upsilon_s^t$ is deterministic and by the filtering property of the conditional expectation we have

$$U_t^{\mathcal{F}}(b) = E[U_t^{\mathcal{G}}(b)|\mathcal{F}_t] = E[U_t^{\mathcal{H}}(b)|\mathcal{F}_t]$$

and in particular

$$(18) \qquad U_0^{\mathcal{F}}(b) = U_0^{\mathcal{G}}(b) = U_0^{\mathcal{H}}(b) = 0.$$

However, as we shall show in subsequent computations, equation (18) does not hold when $\gamma \neq 0$ and, in particular, information neutrality fails to hold in that case.

More explicitly, when $\gamma \neq 0$, one can use the representation (17) and a simple Girsanov transformation to get

$$(19) \qquad U_t^{\mathcal{F}}(b) = \bigg[\int_t^T \Gamma_t^s\,ds\bigg]W_t^{\mathcal{F}} - \gamma\int_t^T (s-t)\Gamma_t^s\,ds,$$

where we used the notation $\Gamma_t^s = \exp(-\int_t^s \beta(u)\,du)$. Thus, by differentiation, the associated intensity is

$$(20) \qquad V_t^{\mathcal{F}}(b) = \int_t^T \Gamma_t^s\,ds.$$

On the other hand, it is also possible to compute the couple $(U^{\mathcal{G}}(b), V^{\mathcal{G}}(b))$. The exercise is slightly more involved and in order to execute it, let us first define a new probability measure on $\mathcal{G}_T$ by

$$\bigg(\frac{d\widetilde{P}}{dP}\bigg)_{\mathcal{G}_T} = \exp\bigg(-\frac{\gamma^2}{2}T - \gamma W_T^{\mathcal{G}}\bigg).$$

By Girsanov's theorem, the process
$$\widetilde{W}_t^{\mathcal{G}} = W_t^{\mathcal{G}} + \gamma t$$
is a $(\widetilde{P}, \mathcal{G})$-Brownian motion and reexpressing the representation (17) under the probability $\widetilde{P}$ and the filtration $\mathcal{G}$ gives
$$\begin{aligned} U_t^{\mathcal{G}}(b) &= \widetilde{E}\Big[\int_t^T \Gamma_t^s \int_0^s \operatorname{sgn}(W_u^{\mathcal{G}})\, dW_u^{\mathcal{G}}\, ds \Big| \mathcal{G}_t\Big] \\ &= \Big[\int_t^T \Gamma_t^s\, ds\Big] W_t^{\mathcal{F}} + \widetilde{E}\Big[\int_t^T \Gamma_t^s \int_t^s \operatorname{sgn}(W_u^{\mathcal{G}})\, dW_u^{\mathcal{G}}\, ds \Big| \mathcal{G}_t\Big], \end{aligned}$$
where $\widetilde{E}$ is the expectation under $\widetilde{P}$. Substituting $W^{\mathcal{G}}$ with $\widetilde{W}^{\mathcal{G}}$ in the above expression and eliminating the stochastic integrals (which does not contribute to the expectation) gives
$$U_t^{\mathcal{G}}(b) = \Big[\int_t^T \Gamma_t^s\, ds\Big] W_t^{\mathcal{F}} - \gamma \int_t^T \Gamma_t^s \int_t^s \widetilde{E}[\operatorname{sgn}(\widetilde{W}_u^{\mathcal{G}} - \gamma u) | \mathcal{G}_t]\, du\, ds.$$

Now, one can remark that
$$\begin{aligned} \widetilde{E}[\operatorname{sgn}(\widetilde{W}_u^{\mathcal{G}} - \gamma u) | \mathcal{G}_t] &= \widetilde{E}[\operatorname{sgn}((\widetilde{W}_u^{\mathcal{G}} - \widetilde{W}_t^{\mathcal{G}}) + \widetilde{W}_t^{\mathcal{G}} - \gamma u) | \mathcal{G}_t] \\ &= \widetilde{E}[\operatorname{sgn}(\sqrt{u-t}\, G + W_t^{\mathcal{G}} - \gamma(u-t))], \end{aligned}$$
where $G$ is a standard Gaussian variable with zero mean and unit variance under the probability $\widetilde{P}$. Consequently, expressing the above quantity in terms of the cumulative $\Phi$ of the standard Gaussian variable gives
$$\widetilde{E}[\operatorname{sgn}(\widetilde{W}_u^{\mathcal{G}} - \gamma u) | \mathcal{G}_t] = 1 - 2\Phi\Big(-\frac{W_t^{\mathcal{G}} - \gamma(u-t)}{\sqrt{u-t}}\Big)$$
and thus
$$\begin{aligned} U_t^{\mathcal{G}}(b) &= \Big[\int_t^T \Gamma_t^s\, ds\Big] W_t^{\mathcal{F}} - \gamma \int_t^T (s-t) \Gamma_t^s\, ds \\ &\quad + 2\gamma \int_t^T \Big(\int_t^s \Phi\Big(-\frac{W_t^{\mathcal{G}} - \gamma(u-t)}{\sqrt{u-t}}\Big)\, du\Big) \Gamma_t^s\, ds. \end{aligned} \tag{21}$$

Differentiating the above expression, and taking only the martingale part gives the intensity
$$\begin{aligned} V_t^{\mathcal{G}}(b) &= \Big[\int_t^T \Gamma_t^s\, ds\Big] \operatorname{sgn}(W_t^{\mathcal{G}}) \\ &\quad - \gamma \int_t^T \Big(\int_t^s \sqrt{\frac{2}{\pi(u-t)}} \exp\Big(-\frac{1}{2}\frac{(W_t^{\mathcal{G}} - \gamma(u-t))^2}{u-t}\Big)\, du\Big) \Gamma_t^s\, ds. \end{aligned}$$

Thus, it becomes clear from (19) and (21) that

$$(22) \quad U_t^{\mathcal{G}}(b) - U_t^{\mathcal{F}}(b) = 2\gamma \int_t^T \left( \int_t^s \Phi\left( -\frac{W_t^{\mathcal{G}} - \gamma(u-t)}{\sqrt{u-t}} \right) du \right) \Gamma_t^s \, ds \neq 0,$$

$P \otimes dt$ a.s. and hence the information neutrality does not hold.

Finally, for the GSDUs under $\mathcal{H}$, one can use some similar Girsanov transformations and get

$$(23) \quad \begin{aligned} U_t^{\mathcal{H}}(b) &= \int_t^T W_s^{\mathcal{F}} \Gamma_t^s \, ds, \\ V_t^{\mathcal{H}}(b) &= 0, \end{aligned}$$

if $t \geq T/2$ and

$$(24) \quad \begin{aligned} U_t^{\mathcal{H}}(b) &= \int_t^{2t} W_s^{\mathcal{F}} \Gamma_t^s \, ds + W_{2t}^{\mathcal{F}} \int_{2t}^T \Gamma_t^s \, ds - \gamma\sqrt{2} \int_{2t}^T (s/2 - t)\Gamma_t^s \, ds, \\ V_t^{\mathcal{H}}(b) &= \sqrt{2} \int_{2t}^T \Gamma_t^s \, ds, \end{aligned}$$

if $t < T/2$.

Thus, we see that

$$(25) \quad U_0^{\mathcal{H}}(b) - U_0^{\mathcal{F}}(b) = -\gamma(1 - \tfrac{\sqrt{2}}{2}) \left( \int_0^T s\Gamma_t^s \, ds \right),$$

and in particular the identity (18) does not hold and hence the information neutrality does not hold again. Note that when $\gamma > 0$, it is clear that $U_0^{\mathcal{H}}(b) < U_0^{\mathcal{F}}(b)$ and we interpret this inequality as a form of aversion to information (the consumer prefers to have access only to the coarser filtration $\mathcal{F}_{(\cdot)}$). However, this inequality is not true for every consumption plan and it can be proved that we have the opposite inequality $U_0^{\mathcal{H}}(b') > U_0^{\mathcal{F}}(b')$ for the consumption plan $b'$ defined by $b'_t = \exp(-W_t^{\mathcal{F}})$.

4.2. *A nonlinear intertemporal aggregator.* Let us consider the Chen and Epstein (2002) GSDU intertemporal aggregator given by

$$(26) \quad f(s, c, y, z) = \log(c) - \beta(s)y - k|z|,$$

for $k \in \mathbb{R}^+$ and for a deterministic integrable function $\beta(\cdot)$.

By Proposition 2, we know that for each filtration $\mathcal{A} \in \{\mathcal{F}, \mathcal{G}\}$, the $\mathcal{A}$-GSDU associated with the intertemporal aggregator (26) maybe represented for any $c \in \mathcal{H}^2(\mathcal{A}, \mathbb{R})$ as

$$(27) \quad U_t^{\mathcal{A}}(c) = \operatorname*{ess\,inf}_{\theta \in \Theta_{\mathcal{A}}} U_t^{\mathcal{A},\theta}(c), \qquad P \otimes dt \text{ a.s.},$$

where $\Theta_{\mathcal{A}} = \{\theta \in \mathcal{H}^2(\mathcal{A}, \mathbb{R}) : |\theta_t| \leq k, \ P \otimes dt \text{ a.s.}\}$ and where the process $U^{\mathcal{A},\theta}(c)$ is the GSDU defined by

$$U_t^{\mathcal{A},\theta}(c) = \int_t^T (\log(c_s) - \beta(s) U_s^{\mathcal{A},\theta}(c) - \theta_s V_s^{\mathcal{A},\theta}(c)) \, ds - \int_t^T V_s^{\mathcal{A},\theta}(c) \, dW_s^{\mathcal{A}},$$

for each $\theta \in \Theta_{\mathcal{A}}$. By Proposition 1, we have the representation

(28) $$U_t^{\mathcal{A},\theta}(c) = E\left[\int_t^T \Upsilon_t^s \log(c_s) \, ds \Big| \mathcal{A}_t\right],$$

where $\Upsilon_t^s$ is the adjoint process defined for $s \geq t$ by the forward SDE

$$d\Upsilon_t^s = -\Upsilon_t^s [\beta_s \, ds + \theta_s \cdot dW_s^{\mathcal{A}}], \qquad \Upsilon_t^t = 1.$$

In order to compute $U_t^{\mathcal{F}}(b)$ where we recall that the consumption plan $b$ is defined in (15), we first define for each $\theta \in \Theta_{\mathcal{F}}$ the probability measure $P^\theta$ by its Radon–Nikodym derivative with respect to $P$ on $\mathcal{F}_T$,

$$\frac{dP^\theta}{dP} = \exp\left\{-\int_0^T \theta_s \, dW_s^{\mathcal{F}} - \frac{1}{2} \int_0^T \theta_s^2 \, ds\right\}.$$

By Girsanov's theorem, the representation (28) becomes

$$U_t^{\mathcal{F},\theta}(b) = E_\theta\left[\int_t^T \Gamma_t^s W_s^{\mathcal{F}} \, ds \Big| \mathcal{F}_t\right]$$
$$= \left[\int_t^T \Gamma_t^s \, ds\right] W_t^{\mathcal{F}} - E_\theta\left[\int_t^T \left(\Gamma_t^s \int_t^s \theta_v \, dv\right) ds \Big| \mathcal{F}_t\right],$$

where $E_\theta$ is the expectation under the probability $P^\theta$. Therefore, the essential infimum of (27) is attained by $\theta = k$, $P \otimes dt$ a.s., and hence

$$U_t^{\mathcal{F}}(b) = \left[\int_t^T \Gamma_t^s \, ds\right] W_t^{\mathcal{F}} - k \int_t^T (s-t) \Gamma_t^s \, ds,$$

and consequently the associated intensity is given by

$$V_t^{\mathcal{F}}(b) = \int_t^T \Gamma_t^s \, ds.$$

In order to compute $U_t^{\mathcal{G}}(b)$, one can write the BSDE satisfied by $U_t^{\mathcal{F}}(b)$, and translate it under $\mathcal{G}$ and get a uniqueness argument

$$U_t^{\mathcal{G}}(b) = \left[\int_t^T \Gamma_t^s \, ds\right] W_t^{\mathcal{F}} - k \int_t^T (s-t) \Gamma_t^s \, ds \equiv U_t^{\mathcal{F}}(b)$$

and

$$V_t^{\mathcal{G}}(b) = \left[\int_t^T \Gamma_t^s \, ds\right] \text{sgn}(W_t^{\mathcal{G}}).$$

For the utility under the filtration $\mathcal{H}$, a uniqueness argument [note that $V^{\mathcal{H}}$ has an invariant sign in equations (23) and (24)] allows us to conclude $U^{\mathcal{H}}$ and $V^{\mathcal{H}}$ are given by equations (23) and (24) with the replacement of $\gamma$ by $k$.

4.3. *A quadratic intertemporal aggregator.* Let us consider the SDU intertemporal aggregator given by

(29) $$f(s,c,y,z) = \log(c) - \frac{\gamma}{2}z^2,$$

for $\gamma \in \mathbb{R}$. Although the intertemporal aggregator (29) is not Lipshitz, using Itô's rule, it can be shown easily that

$$U_t^{\mathcal{A}}(c) = -\frac{1}{\gamma}\log\bigg(E\bigg[\exp\bigg(-\gamma\int_t^T \log(c_s)\,ds\bigg)\bigg|\mathcal{A}_t\bigg]\bigg),$$

for any $c \in c \in \mathcal{H}^2(\mathcal{F}, \mathbb{R})$ and for $\mathcal{A} = \mathcal{F}, \mathcal{G}$ and $\mathcal{H}$. In particular, existence and uniqueness of the associated BSDE hold when the above expectation is finite. Furthermore, some straightforward computations give

$$U_t^{\mathcal{F}}(b) = U_t^{\mathcal{G}}(b) = (T-t)W_t^{\mathcal{F}} - \frac{\gamma}{6}(T-t)^3,$$

$$U_t^{\mathcal{H}}(b) = \int_t^{2t} W_s^{\mathcal{F}}\,ds + (T-2t)W_{2t}^{\mathcal{F}} - \frac{\gamma}{6}(T-2t)^3 \qquad \text{for } t < T/2$$

and

$$V_t^{\mathcal{F}}(b) = (T-t),$$
$$V_t^{\mathcal{G}}(b) = (T-t)\operatorname{sgn}(W_t^{\mathcal{G}}),$$
$$V_t^{\mathcal{H}}(b) = \sqrt{2}(T-t).$$

4.4. *Discussion.* While the computations of this section offer a modest contribution from a theoretical perspective, they have the merit of illustrating in a concrete way how a BSDE depends on its filtration.

For our objective of characterizing information neutrality, it is worthwhile to see what we can learn from this example. First, the linear GSDU computations of Section 4.1 show that the dependency in $z$ does not allow information neutrality to hold (for the three types of information heterogeneity under consideration).

When the utility is linear but independent from $z$, the GSDUs under $\mathcal{F}$ and $\mathcal{G}$ coincide but are different from the GSDU under $\mathcal{H}$. This fact suggests that the information heterogeneity $\mathcal{F}$ versus $\mathcal{G}$ has a special feature. This special feature seems to be confirmed by the nonlinear intertemporal aggregators of Sections 4.2 and 4.3. In both cases, the GSDUs under $\mathcal{F}$ and $\mathcal{G}$ coincide but are different from the GSDU under $\mathcal{H}$.

Although, these statements have no theoretical value since they are only valid for one particular consumption plan (b), the next theoretical work will establish that the information heterogeneity of the type $\mathcal{F}$ versus $\mathcal{G}$ is the unique type of information heterogeneity which allow information neutrality to hold for a class of intertemporal aggregators.

**5. Characterization of information neutrality.** In this section we fix a couple of filtrations $\mathcal{F}_{(\cdot)} \subsetneq \mathcal{G}_{(\cdot)}$. First, let us state a general necessary condition of information neutrality.

LEMMA 1. *If a BSDE* (3) [*resp. a GSDU* (5)] *exhibits information neutrality then for any* $\xi \in L^2(\mathcal{F}_T)$ [*resp. for any* $c \in \mathcal{H}^2(\mathcal{F}, \mathbb{R})$] *we have*

$$\|Z_t^{\mathcal{F}}(\xi)\| = \|Z_t^{\mathcal{G}}(\xi)\| \qquad [resp. \ \|V_t^{\mathcal{F}}(c)\| = \|V_t^{\mathcal{G}}(c)\|], \qquad P \otimes dt \ a.s.$$

PROOF. It follows from (7) that for each $\xi \in L^2(\mathcal{F}_T)$,

$$\lim_{|t_{i+1}-t_i|\to 0} \sum_i |Y_{t_{i+1}}^{\mathcal{F}}(\xi) - Y_{t_i}^{\mathcal{F}}(\xi)|^2 = \lim_{|t_{i+1}-t_i|\to 0} \sum_i |Y_{t_{i+1}}^{\mathcal{G}}(\xi) - Y_{t_i}^{\mathcal{G}}(\xi)|^2$$

and therefore,

$$\int_0^t \|Z_s^{\mathcal{F}}(\xi)\|^2 \, ds = \int_0^t \|Z_s^{\mathcal{G}}(\xi)\|^2 \, ds.$$

Consequently,

$$\|Z_t^{\mathcal{F}}(\xi)\| = \|Z_t^{\mathcal{G}}(\xi)\|, \qquad P \otimes dt \ \text{a.s.},$$

and the proof is similar for the GSDU case. □

5.1. *The BSDE problem.* In this section, and in view of Lemma 1, we use the following assumption on the driver of the BSDE (3):

(H1) The BSDE driver has the form

$$h(s, \omega, y, z) = \tilde{h}(s, \omega, y, \|z\|).$$

THEOREM 1. *Under assumption* (H1) *the following statements are equivalent*:

(a) *The BSDE* (3) *exhibits strong information neutrality.*
(b) *There exists a process* $M : [0,T] \times \Omega \to \mathbb{R}^{n \times n}$ *in the set* $\mathcal{P}^{\mathcal{G}}$ *such that* $M'M = Id_n$, $dt \otimes dP$ *a.s. and*

$$(30) \qquad W_t^{\mathcal{F}} = \int_0^t M_s \, dW_s^{\mathcal{G}}, \qquad P \otimes dt \ a.s.$$

PROOF. (b) $\Rightarrow$ (a) For any $\xi \in L^2(\mathcal{F}_T)$, assumption (H1) in conjunction with (30) implies that $(Y_t^{\mathcal{F}}(\xi), M_t Z_t^{\mathcal{F}}(\xi))_{0 \le t \le T}$ solves the BSDE (3), and by uniqueness we get $Y_t^{\mathcal{F}}(\xi) = Y_t^{\mathcal{G}}(\xi)$.

(b) $\Rightarrow$ (a) For any $\xi \in L^2(\mathcal{F}_T)$, and under assumption (H1), substracting (3) and (4) gives, by Lemma 1,

$$(31) \qquad \int_0^t Z_s^{\mathcal{F}}(\xi) \cdot dW_s^{\mathcal{F}} = \int_0^t Z_s^{\mathcal{G}}(\xi) \cdot dW_s^{\mathcal{G}}, \qquad P \otimes dt \ \text{a.s.}$$

In order to compute explicitly $Z^{\mathcal{F}}$ for some particular $\xi$, let us now introduce the $n$-dimensional process $(X_t = (X_t, \ldots, X_t))_{0 \le t \le T}$ that solves the stochastic differential equation

$$dX_t = -g(t, \omega, X_t)\,dt + dW_t^{\mathcal{F}},$$
$$X_0 = x \in \mathbb{R}^n,$$

where $g$ maps $[0,T] \times \Omega \times \mathbb{R}^n$ onto $\mathbb{R}^n$ and is defined by

$$g(t, \omega, X) = (\tilde{h}(t, \omega, X_1, 1); \tilde{h}(t, \omega, X_2, 1); \ldots; \tilde{h}(t, \omega, X_n, 1)).$$

By construction, it is clear that, for each $k = 1, \ldots, n$; $Y_t^{\mathcal{F}}(X_T^k) = X_t^k$ and $Z_t^{\mathcal{F}}(X_T^k) = \delta_k$ where $\delta_k$ is a vector of $\mathbb{R}^n$ defined by $\delta_{kk'} = 0$ if $k \ne k'$ and $\delta_{kk} = 1$. Therefore, for each $k, i = 1, \ldots, n$, letting

$$M_t^{k,i} := Z_t^{i,\mathcal{G}}(X_T^k),$$

we get, from (31),

$${}^k W_t^{\mathcal{F}} = \int_0^t \sum_{i=1}^n M_s^{k,i}\,d\,{}^i W_s^{\mathcal{G}}, \qquad P \otimes dt \text{ a.s.},$$

and thus for each $k, k' = 1, \ldots, n$,

$$\delta_k' \delta_{k'}\,dt = d\langle {}^k W^{\mathcal{F}}, {}^{k'} W^{\mathcal{F}} \rangle_t = \sum_{i=1}^n M_t^{k,i} M_t^{k',i}\,dt = (M_t' M_t)_{k,k'}\,dt,$$

which completes the proof. $\square$

To be concrete, in the scalar Brownian motion case, the filtration $\mathcal{F}_t = \sigma(|W_s^{\mathcal{G}}|, 0 \le s \le t)$ is generated by the Brownian motion $W_t^{\mathcal{F}} := \int_0^t \text{sgn}(W_s^{\mathcal{G}})\,dW_s^{\mathcal{G}}$ as explained in Section 3.1 and this is an example of a situation where condition (b) of Theorem 1 holds. More generally, Theorem 1 illustrates that when the link (30) exists between two filtrations $\mathcal{F}$ and $\mathcal{G}$, there is *no utility cost* for the information loss due to accessing to $\mathcal{F}$ rather than $\mathcal{G}$ under assumption (H1).

REMARK 1. Note that under assumption (H1), it is easily seen from (30) that information neutrality is also equivalent to the fact that $W^{\mathcal{F}}$ is a $\mathcal{G}$-Brownian motion. Thus Theorem 1 provides a possible interpretation of this condition in terms of utility cost of information.

REMARK 2. When $h = 0$, the BSDE (3) is the (linear) conditional expectation and the strong information neutrality becomes

(32) $$E[\xi | \mathcal{F}_t] = E[\xi | \mathcal{G}_t], \qquad dP \otimes dt \text{ a.s.},$$

for all $\xi \in L^2(\mathcal{F}_T)$. In particular, the above is true for $\xi = W_T^{\mathcal{F}}$, and from the Lévy criteria, we deduce that, in fact, $W^{\mathcal{F}}$ is a $\mathcal{G}$-Brownian motion and thus (30) is satisfied. Thus Theorem 1 may be interpreted as a generalization of property (32) to the BSDE generated by the driver $h$ under assumption (H1). (The BSDEs are denominated sometimes nonlinear expectations.)

5.2. *The GSDU case.* In this section, and in view of Lemma 1, we use the following assumption on the GSDU intertemporal aggregator:

(H2) The intertemporal aggregator has the form
$$f(s,c,y,z) = \tilde{f}(s,c,y,\|z\|).$$

Furthermore, for technical reasons we shall use the following assumption on the intertemporal aggregator:

(H3) The intertemporal aggregator $f$ is continuously differentiable with respect to $c, y, z$ with first derivative being bounded by some constant $L > 0$ and satisfies for all $(t, c, y, z) \in [0, T] \times \mathbb{R} \times \mathbb{R} \times \mathbb{R}^n$
$$|f(t,c,0,z)| \leq C, \qquad |\partial_c f(t,c,y,z)| \geq k > 0,$$
for some constants $C$ and $k$. Moreover, the derivatives $\partial_c f, \partial_y f, \partial_z f$ are uniformly Lipschitz with respect to each of the variables $c, y, z$ with a Lipschitz constant $M > 0$.

THEOREM 2. *Under assumptions* (H2) *and* (H3), *the following statements are equivalent*:

(a) *The GSDU* (5) *exhibits information neutrality.*
(b) *Equality* (30) *holds for some process* $M : [0, T] \times \Omega \to \mathbb{R}^{n \times n}$ *in the set* $\mathcal{P}^{\mathcal{G}}$ *such that* $M'M = Id_n$, $dt \otimes dP$ *a.s.*

PROOF. (b) $\Rightarrow$ (a) For any $c \in \mathcal{H}^2(\mathcal{F}, \mathbb{R})$, assumption (H2) in conjunction with (30) implies that $(U_t^{\mathcal{F}}(c), M_t'V_t^{\mathcal{F}}(c))_{0 \leq t \leq T}$ solves the BSDE (5), and by uniqueness we get $U^{\mathcal{F}}(c) = U^{\mathcal{G}}(c)$.

(a) $\Rightarrow$ (b) A similar approach of Theorem 1 leads under assumption (H2) to the identity

$$(33) \qquad \int_0^t V_s^{\mathcal{F}}(c)\, dW_s^{\mathcal{F}} = \int_0^t V_s^{\mathcal{G}}(c)\, dW_s^{\mathcal{G}}, \qquad P \otimes dt \text{ a.s.}$$

for each $c \in \mathcal{H}^2(\mathcal{F}, \mathbb{R})$.

In particular, for $k = 1, \ldots, n$, let us consider the utility $U_t^{\mathcal{F}}(^kW^{\mathcal{F}})$ associated to the consumption processes $c_t = {}^kW_t^{\mathcal{F}}$ that coincides, by uniqueness of the BSDE (5), with the solution of the scalar BSDE

$$^kN_t = \int_t^T \tilde{f}(s, {}^kW_s^{\mathcal{F}}, {}^kN_s, |{}^k\zeta_s|)\, ds - \int_t^T {}^k\zeta_s\, d^kW_s^{\mathcal{F}}.$$

Hence, $V_t^{\mathcal{F}}(^k W^{\mathcal{F}}) = {}^k\zeta_t \delta_k$, where we recall that $\delta_k$ is a vector of $\mathbb{R}^n$ defined by $\delta_{kk'} = 0$ if $k \neq k'$ and $\delta_{kk} = 1$. Now, assuming that

(34) $$|{}^k\zeta_t| > 0, \qquad P \otimes dt \text{ a.s. for } k = 1, \ldots, n,$$

and applying (33) to the consumption processes $c = {}^k W^{\mathcal{F}}$ for $k = 1, \ldots, n$ gives

$${}^k W_t^{\mathcal{F}} = \int_0^t \sum_{i=1}^n M_s^{k,i} d^i W_s^{\mathcal{G}}, \qquad P \otimes dt \text{ a.s.},$$

where

$$M_t^{k,i} := \frac{V_t^{i,\mathcal{G}}({}^k W^{\mathcal{F}})}{{}^k\zeta_t}$$

and following the same argument of Theorem 1 we are done. The following lemma shows that under assumption (H3), inequality (34) is satisfied. □

LEMMA 2. *Let $(B_t, 0 \leq t \leq T)$ be a standard unidimensional Brownian motion [say under $(P, \mathcal{F})$] and consider the BSDE*

$$N_t^0 = \int_t^T g(s, B_s, N_s^0, \zeta_s^0) \, ds - \int_t^T \zeta_s^0 \, dB_s,$$

*where $g$ is defined for each $(t, c, y, z) \in [0, T] \times \mathbb{R} \times \mathbb{R} \times \mathbb{R}$ by*

$$g(t, c, y, z) = f(t, c, y, z\delta_1)$$

*and where $f$ is an intertemporal aggregator satisfying assumption* (H3).
*Then $|\zeta_t^0| > 0, P \otimes dt$ a.s.*

PROOF. First, consider the finality of BSDEs parametrized by $x \geq 0$ and defined by

(35) $$N_t^x = \int_t^T g(s, x + B_s, N_s^x, \zeta_s^x) \, ds - \int_t^T \zeta_s^x \, dB_s.$$

By a result of Ma, Protter and Yong [(1994), Lemma 3.2; see also Pardoux and Peng (1992)] and for each $x \geq 0$ the solution of (35) satisfies

$$(N_t^x, \zeta_t^x) = (\theta(t, x + B_t), \partial_x \theta(t, x + B_t)),$$

where $\theta : [0, T] \times \mathbb{R} \to \mathbb{R}$ is the unique bounded classical solution of the quasi-linear parabolic equation

$$\partial_t \theta(t, x) + \tfrac{1}{2} \partial_{xx} \theta(t, x) + g(t, x, \theta(t, x), \partial_x \theta(t, x)) = 0, \qquad (t, x) \in (0, T) \times \mathbb{R},$$
$$\theta(T, x) = 0.$$

Furthermore, by the BSDE a priori estimates [see Pardoux and Peng (1990) and El Karoui, Peng and Quenez (1997)] and the boundness of $\partial_c g$ [assumption (H3)] we have the following bounding argument:

$$
\begin{aligned}
E&\left(\sup_{t\le T}|N_t^x - N_t^0|^2\right) + E\left(\int_0^T dt\,(\zeta_t^x - \zeta_t^0)^2\right) \\
(36)\quad &\le CE\left(\int_0^T dt,\ (g(t, x+B_t, N_t^0, \zeta_t^0) - g(t, B_t, N_t^0, \zeta_t^0))^2\right) \\
&\le CTL^2 x^2,
\end{aligned}
$$

for some positive constant $C$.

Second, consider the linear BSDE

$$
\begin{aligned}
(37)\quad Q_t = \int_t^T &(\partial_c g(s, B_s, N_s^0, \zeta_s^0) \\
&+ Q_s\,\partial_y g(s, B_s, N_s^0, \zeta_s^0) + \Lambda_s\,\partial_z g(s, B_s, N_s^0, \zeta_s^0))\,ds \\
&- \int_t^T \Lambda_s\,dB_s
\end{aligned}
$$

which, according to Proposition 1, gives the following solution:

$$
(38)\quad Q_t = E\left[\int_t^T ds\,\partial_c g(s, B_s, N_s^0, \zeta_s^0)\exp\left(\int_t^s du\,\partial_y g(u, B_u, N_u^0, \zeta_u^0)\right)\Upsilon_t^s\Big|\mathcal{F}_t\right],
$$

where

$$
\Upsilon_t^s = \exp\left(-\tfrac{1}{2}\int_t^s du\,\partial_z g(u, B_u, N_u^0, \zeta_u^0)^2 + \int_t^s dB_u\,\partial_z g(u, B_u, N_u^0, \zeta_u^0)\right).
$$

Notice that from the boundness of the first derivatives of $g$, equation (38) shows that $Q$ satisfies

$$
(39)\quad ke^{-MT}(T - t) \le |Q_t| \le Me^{MT}T
$$

and, therefore, in order to prove the lemma, we will show that

$$
(40)\quad \zeta_t^0 = Q_t, \qquad P \otimes dt \text{ a.s.}
$$

More precisely, since $\zeta_t^0 = \partial_x \theta(t, B_t)$ and by the definition of a derivative, we need only prove

$$
(41)\quad \lim_{x\downarrow 0}\Delta N_t^x = Q_t, \qquad P \otimes dt \text{ a.s.},
$$

where $\Delta N_t^x := (N_t^x - N_t^0)/x$ for $x > 0$.

To that end, let us define for each $x > 0$ the process $\Delta\zeta_t^x := (\zeta_t^x - \zeta_t^0)/x$ and, following a linearization technique of El Karoui, Peng and Quenez

(1997) we interpret the couple $(\Delta N_t^x, \Delta \zeta_t^x)$ as the solution of the linear BSDE

$$\Delta N_t^x = \int_t^T (\varphi_s^x + A_s^x \Delta N_s^x + B_s^x \Delta \zeta_s^x) \, ds - \int_t^T \Delta \zeta_s^x \, dB_s,$$

where

$$A_t^x := \int_0^1 d\lambda \, \partial_y g(t, x + B_t, N_t^0 + \lambda(N_t^x - N_t^0), \zeta_t^x),$$

$$B_t^x := \int_0^1 d\lambda \, \partial_z g(t, x + B_t, N_t^0, \zeta_t^0 + \lambda(\zeta_t^x - \zeta_t^0)),$$

$$\varphi_t^x := (g(t, x + B_t, N_t^0, \zeta_t^0) - g(t, B_t, N_t^0, \zeta_t^0))/x.$$

Using the BSDE a priori estimates [see Pardoux and Peng (1990) and El Karoui, Peng and Quenez (1997)] and the inequality $(a + b + c)^2 \leq 4(a^2 + b^2 + c^2)$ we get

$$
\begin{aligned}
E\bigg(&\sup_{t \leq T} |\Delta N_t^x - Q_t|^2\bigg) \\
&\leq CE\bigg(\int_0^T dt \, (A_t^x - \partial_y g(t, B_t, N_t^0, \zeta_t^0))^2 Q_t^2\bigg) \\
&\quad + CE\bigg(\int_0^T dt \, (B_t^x - \partial_z g(t, B_t, N_t^0, \zeta_t^0))^2 \Lambda_t^2\bigg) \\
&\quad + CE\bigg(\int_0^T dt \, (\varphi_t^x - \partial_c g(t, B_t, N_t^0, \zeta_t^0))^2\bigg),
\end{aligned}
\tag{42}
$$

for some positive constant $C$.

By Rolle's theorem, we have $\varphi_t^x = \partial_c g(t, \eta_t^x, N_t^0, \zeta_t^0)$ for some $\eta_t^x \in (B_t, B_t + x)$ and since $\partial_c g$ is uniformly Lipschitz with respect to $c$, we have the following bound for the third term on the right-hand side of (42):

$$CE\bigg(\int_0^T dt \, (\varphi_t^x - \partial_c g(t, B_t, N_t^0, \zeta_t^0))^2\bigg) \leq CTM^2 x^2. \tag{43}$$

To analyze the convergence of the first term on the right-hand side of (42) we need the following bounding arguments:

$$
\begin{aligned}
CE\bigg(\int_0^T dt \bigg(\int_0^1 d\lambda \, \partial_y g(t, x + B_t, N_t^0 + \lambda(N_t^x - N_t^0), \zeta_t^x) \\
- \partial_y g(t, B_t, N_t^0, \zeta_t^0)\bigg)^2 Q_t^2\bigg) \\
\leq KE\bigg(\int_0^T dt \bigg(\int_0^1 d\lambda \, \partial_y g(t, x + B_t, N_t^0 + \lambda(N_t^x - N_t^0), \zeta_t^x)
\end{aligned}
$$

$$-\partial_y g(t, B_t, N_t^0, \zeta_t^0)\Big)^2\Big)$$

$$\leq 4KE\bigg(\int_0^T dt \bigg(\int_0^1 d\lambda\, \partial_y g(t, x+B_t, N_t^0 + \lambda(N_t^x - N_t^0), \zeta_t^x)$$

$$- \partial_y g(t, x+B_t, N_t^0, \zeta_t^x)\Big)^2\bigg)$$

$$+ 4KE\bigg(\int_0^T dt \bigg(\int_0^1 d\lambda\, \partial_y g(t, x+B_t, N_t^0, \zeta_t^x)$$

$$- \partial_y g(t, x+B_t, N_t^0, \zeta_t^0)\Big)^2\bigg)$$

$$+ 4KE\bigg(\int_0^T dt \bigg(\int_0^1 d\lambda\, \partial_y g(t, x+B_t, N_t^0, \zeta_t^0)$$

$$- \partial_y g(t, B_t, N_t^0, \zeta_t^0)\Big)^2\bigg)$$

$$\leq 4KM^2\bigg(E\bigg(\int_0^T dt\, (N_t^x - N_t^0)^2/3\bigg) + E\bigg(\int_0^T dt\, (\zeta_t^x - \zeta_t^0)^2\bigg) + x^2\bigg)$$

$$\leq 4KM^2(CT^2L^2/3 + CTL^2 + 1)x^2,$$

where $K = CM^2 e^{2MT} T^2$ and where we used (39) to obtain the first inequality, the uniform Lipschitz property of $\partial_y f$ to obtain third inequality and we used the a priori estimates in (36) to obtain the last inequality. By (43) and the above, we therefore have, on taking limits,

(44)
$$\lim_{x \downarrow 0} CE\bigg(\int_0^T dt\, (A_t^x - \partial_y g(t, B_t, N_t^0, \zeta_t^0))^2 Q_t^2\bigg)$$
$$+ CE\bigg(\int_0^T dt\, (\varphi_t^x - \partial_c g(t, B_t, N_t^0, \zeta_t^0))^2\bigg) = 0.$$

In order to tackle the convergence of the second term on the right-hand side of (42) we use similar bound to the above and get

$$CE\bigg(\int_0^T dt \bigg(\int_0^1 d\lambda\, \partial_z g(t, x+B_t, N_t^0, \zeta_t^0 + \lambda(\zeta_t^x - \zeta_t^0))$$

$$- \partial_z g(t, B_t, N_t^0, \zeta_t^0)\Big)^2 \Lambda_t^2\bigg)$$

$$\leq 2CE\bigg(\int_0^T dt \bigg(\int_0^1 d\lambda\, \partial_z g(t, x+B_t, N_t^0, \xi_t^0 + \lambda(\zeta_t^x - \zeta_t^0))$$

$$- \partial_z g(t, x+B_t, N_t^0, \zeta_t^0)\Big)^2 \Lambda_t^2\bigg)$$

$$(45) \quad + 2CE\left(\int_0^T dt \left(\int_0^1 d\lambda\, \partial_z g(t, x+B_t, N_t^0, \zeta_t^0)\right.\right.$$
$$\left.\left. - \partial_z g(t, B_t, N_t^0, \zeta_t^0)\right)^2 \Lambda_t^2\right)$$
$$\leq 2CE\left(\int_0^T dt \left(\int_0^1 d\lambda\, \partial_z g(t, x+B_t, N_t^0, \zeta_t^0 + \lambda(\zeta_t^x - \zeta_t^0))\right.\right.$$
$$\left.\left. - \partial_z g(t, x+B_t, N_t^0, \zeta_t^0)\right)^2 \Lambda_t^2\right)$$
$$+ 2M^2 x^2 E\left(\int_0^T dt\, \Lambda_t^2\right).$$

Finally, to bound the last term of (45) we mimic a technique introduced in El Karoui, Peng and Quenez (1997) as follows:

$$E\left(\int_0^T dt \left(\int_0^1 d\lambda\, \partial_z g(t, x+B_t, N_t^0, \zeta_t^0 + \lambda(\zeta_t^x - \zeta_t^0))\right.\right.$$
$$\left.\left. - \partial_z g(t, x+B_t, N_t^0, \zeta_t^0)\right)^2 \Lambda_t^2\right)$$
$$\leq M^2 x^{1/2} E\left(\int_0^T dt\, \Lambda_t^2\right)$$
$$+ E\left(\int_0^T dt \left(\int_0^1 d\lambda\, \Lambda_t^2 \mathbb{1}_{\{|\zeta_t^x - \zeta_t^0| > x^{1/4}\}}\right.\right.$$
$$\times\, \partial_z g(t, x+B_t, N_t^0, \zeta_t^0 + \lambda(\zeta_t^x - \zeta_t^0))$$
$$\left.\left. - \partial_z g(t, x+B_t, N_t^0, \zeta_t^0)\right)^2\right)$$
$$(46) \quad \leq M^2 x^{1/2} E\left(\int_0^T dt\, \Lambda_t^2\right) + 4L^2 E\left(\int_0^T dt\, \mathbb{1}_{\{|\zeta_t^x - \zeta_t^0| > x^{1/4}\}} \Lambda_t^2\right)$$
$$\leq M^2 x^{1/2} E\left(\int_0^T dt\, \Lambda_t^2\right)$$
$$+ 4L^2 E\left(\int_0^T dt\, \mathbb{1}_{\{|\Lambda_t| > x^{-1/4}\}} \mathbb{1}_{\{|\zeta_t^x - \zeta_t^0| > x^{1/4}\}} \Lambda_t^2\right)$$
$$+ 4L^2 x^{-1/2} E\left(\int_0^T dt\, \mathbb{1}_{\{|\zeta_t^x - \zeta_t^0| > x^{1/4}\}}\right)$$
$$\leq M^2 x^{1/2} E\left(\int_0^T dt\, \Lambda_t^2\right) + 4L^2 E\left(\int_0^T dt\, \mathbb{1}_{\{|\Lambda_t| > x^{-1/4}\}} \Lambda_t^2\right)$$

$$+ 4L^2 x^{-1} E\left(\int_0^T dt |\zeta_t^x - \zeta_t^0|^2\right)$$

$$\leq M^2 x^{1/2} E\left(\int_0^T dt\, \Lambda_t^2\right) + 4L^2 E\left(\int_0^T dt\, \mathbb{1}_{\{|\Lambda_t| > x^{-1/4}\}} \Lambda_t^2\right) + 4CTL^4 x,$$

where we used the uniform Lipschitz and the boundness of $\partial_z g$ to obtain the first and the second inequality, the Markov inequality to obtain fourth inequality and (36) to obtain the last inequality. Next, since $\Lambda_t$ is square integrable by construction, the Lebesgue theorem implies that

$$\lim_{x \downarrow 0} E\left(\int_0^T dt\, \mathbb{1}_{\{|\Lambda_t| > x^{-1/4}\}} \Lambda_t^2\right) = 0,$$

and we conclude that the bound obtained in (46) tends to 0 when $x \downarrow 0$, and thus so does the bound obtained in (45).

The above analysis, in conjunction with (44), implies

$$\lim_{x \downarrow 0} E\left(\sup_{t \leq T} |\Delta N_t^x - Q_t|^2\right) = 0$$

and, in particular, (41) and thus (40) are established and we are done. $\square$

Similarly to the BSDE model, note that under assumptions (H2) and (H3), information neutrality for the GSDU model is equivalent to the fact that $W^\mathcal{F}$ is a $\mathcal{G}$-Brownian motion. Thus Theorem 1 provides a possible interpretation of this condition in terms of utility cost of information in the formal sense of the GSDU model. From that perspective, the results of Theorems 1 and 2 have similar interpretations.

The technical proofs of Theorems 1 and 2 are also similar, except in the part constructing a rich set of $\mathcal{F}$-martingales. In the BSDE model (Theorem 1), we used forward SDEs to construct efficiently an appropriate set of terminal data ($\xi$) that generate $n$ $\mathcal{F}$-martingales which integrands form a basis of $\mathbb{R}^n$. This last property allowed us to invert the time derivative of formula (31) and identify $dW^\mathcal{F}$ in terms of $dW^\mathcal{G}$. In the BSDE model, we did not need extra technical conditions because the assumed Lipschitz conditions on the driver $h$ are sufficient to secure the existence of strong solutions of the forward SDEs.

In the GSDU case, the forward SDEs technique will not work because the terminal data is fixed ($\xi = 0$) and we can only choose the consumption process. In other words, we need to control the intensity $V$ of the GSDU by selecting appropriately the intertemporal aggregators [through the choice of the consumption process ($c$)]. The method that we provide in Lemma 2 is a self-contained proof which relies on the link between BSDEs and quasilinear parabolic differential equations [Pardoux and Peng (1992) and Ma, Protter and Yong ([1994](#))]. This method requires extra technical conditions [assumption (H3)] but it is our beliefs that it is interesting by itself and it also has the merit to rely only on classical results from the BSDE literature.

However, assumption (H3) is intuitively not a necessary condition to obtain our characterization [unlike assumption (H2)] and it is *possible* to weaken it. For instance, using the representation of the intensity as a right limit of the Malliavin derivatives of the utility process ($V_t = \lim_{s \downarrow t} \mathcal{D}_t U_s$), one can use the results of Pardoux and Peng ([1992](#)) and El Karoui, Peng and Quenez ([1997](#)) [see also Ma and Zhang ([2002](#))] to weaken assumption (H3). More generally, our choice of consumption process $c_t = {}^k W_t^{\mathcal{F}}$ to generate a rich set of $\mathcal{F}$-martingales is particular. Arguably it is also possible to weaken the required technical conditions of any method (either the Markovian technique or the Malliavin derivative technique) by a different choice of consumption process.

**6. Conclusion.** We characterized the information neutrality property for a class of BSDEs including GSDUs under the assumption that the driver depends on the intensity $Z$ only through the Euclidian norm $\|Z\|$. Behaviorally, the information neutrality property corresponds to a form of intrinsic indifference to information. We proved that, unless the information reduction is specific, the class of GSDUs exhibits an intrinsic attitude toward information. This intrinsic attitude toward information is in fact inflexibly associated to the risk aversion and the ambiguity aversion concepts and cannot be disentangled from them within the GSDU context. These results invite further analysis. In particular, it would be meaningful to characterize the monotonicity, that is, strict preference (or aversion) for information. The Kreps–Porteus–Skiadas approach only provides sufficient conditions for the monotonicity for a particular class of GSDUs: the SDUs. New techniques must be introduced since we have to manipulate sub-solution of BSDEs. For instance the generating martingale technique will not be useful in this context because we only have inequalities.

**Acknowledgments.** I am grateful to Shiqi Song for the very interesting discussions that we had on this topic and to Francis Hirsch for the special support that he gave me for this work. I thank Monique Jeanblanc and Tan Wang for their helpful comments. I thank the anonymous referees for their helpful insight in a first version of this paper (previously entitled "Information Neutrality in the Stochastic Differential Utility").

## REFERENCES

Anderson, E. W., Hansen, L. P. and Sargent, T. J. (1998). Risk and robustness in equilibrium. Working paper, Dept. Economics, Univ. Chicago.

Artzner, P., Delbaen, F., Eber, J. M. and Heath, D. (1999). Coherent measures of risk. *Math. Finance* **4** 203–228. MR1850791

Artzner, P., Delbaen, F., Eber, J. M., Heath, D. and Ku, H. (2002). Coherent multiperiod risk measurement. Working paper, ETH, Zurich.

Chen, Z. and Epstein, L. (2002). Ambiguity, risk and asset returns in continuous time. *Econometrica* **70** 1403–1443. MR1929974

Chew, S.-H. and Ho, J. L. (1994). Hope: An empirical study of the attitude toward the timing of uncertainty resolution. *Journal of Risk and Uncertainty* **8** 267–288.

Duffie, D. and Epstein, L. (1992). Stochastic differential utility. *Econometrica* **60** 353–394. MR1162620

El Karoui, N., Peng, S. and Quenez, M.-C. (1997). Backward stochastic differential equations in finance. *Math. Finance* **7** 1–71. MR1434407

Epstein, L. (1980). Decision making and the temporal resolution of uncertainty. *Internat. Econom. Rev.* **21** 269–283. MR581954

Epstein, L. and Zin, S. (1989). Substitution, risk aversion and the temporal behavior of consumption and asset returns. *Econometrica* **57** 937–969. MR1006550

Epstein, L. and Zin, S. (1991). Substitution, risk aversion and the temporal behavior of consumption and asset returns: An empirical analysis. *Journal of Political Economy* **99** 263–286.

Grant, S., Kajii, A. and Polack, B. (1998). Intrinsic preference for information. *J. Econom. Theory* **83** 233–259. MR1664058

Hansen, L. P., Sargent, T., Turmuhambetova, G. and Williams, N. (2002). Robustness and uncertainty aversion. Working paper, Dept. Economics, Univ. Chicago.

Karatzas, I. and Shreve, S. E. (1988). *Brownian Motions and Stochastic Calculus*. Springer, Berlin. MR917065

Kobylanski, M. (2000). Backward stochastic differential equations and partial differential equations with quadratic growth. *Ann. Probab.* **28** 558–602. MR1782267

Kreps, D. M. and Porteus, E. L. (1978). Temporal resolution of uncertainty and dynamic choice theory. *Econometrica* **46** 185–200. MR481590

Lazrak, A. and Quenez, M. C. (2003). A generalized stochastic differential utility. *Math. Oper. Res.* **28** 154–180. MR1961272

Ma, J. and Yong, J. (1999). *Forward-Backward Stochastic Differential Equations and Their Applications. Lecture Notes in Math.* **1702**. Springer, Berlin. MR1704232

Ma, J. and Zhang, J. (2002). Representation theorems for backward stochastic differential equations. *Ann. Appl. Probab.* **12** 1390–1418. MR1936598

Ma, J., Protter, P. and Yong, J. M. (1994). Solving forward-backward stochastic differential equations explicitly-a four step scheme. *Probab. Theory Related Fields* **98** 339–359. MR1262970

Liptser, R. S. and Shiryayev, A. N. (1977). *Statistics of Random Processes* **I**. Springer, Berlin. MR474486

Pardoux, E. and Peng, S. (1990). Adapted solution of a backward stochastic differential equation. *Systems Control Lett.* **14** 55–61. MR1037747

Pardoux, E. and Peng, S. (1992). Backward stochastic differential equation and quasilinear parabolic partial differential equations. *Stochastic Partial Differential Equations and Their Applications. Lecture Notes in Control and Inform. Sci.* **176** 200–217. Springer, Berlin. MR1176785

Revuz, D. and Yor, M. (1999). *Continuous Martingales and Brownian Motion*, 3rd ed. Springer, Berlin. MR1725357

Riedel, F. (2002). Dynamic coherent risk measures. Working paper, Bonn Univ. MR1892189

Shroder, M. and Skiadas, C. (1999). Optimal consumption and portfolio selection with stochastic differential utility. *J. Econom. Theory* **89** 68–126. MR1724377

Skiadas, C. (1998). Recursive utility and preferences for information. *Econom. Theory* **12** 293–312. MR1658816

Skiadas, C. (2003). Robust control and recursive utility. *Finance Stoch.* **7** 475–489. MR2014246

Uppal, R. and Wang, T. (2003). Model misspecification and underdiversification. *J. Finance* **58** 2465–2486.

Wang, T. (2000). A class of dynamic risk measures. Preprint, Univ. British Columbia.

University of British Columbia
Sauder School of Business
2053 Main Mall
Vancouver, British Columbia
Canada V6T 1Z2
e-mail: lazrak@sauder.ubc.ca